# Estimating of the number of natural solutions of homogeneous algebraic Diophantine diagonal equations with integer coefficients

VICTOR VOLFSON

ABSTRACT. Author developed a method in the paper, which, unlike the circle method of Hardy and Littlewood (CM), allows you to perform a lower estimate for the number of natural (integer) solutions of algebraic Diophantine equation with integer coefficients. It was found the lower estimate of the number of natural solutions to various types of homogeneous algebraic Diophantine equations with integer coefficients diagonal form with any number of variables using this method. Author obtained upper bound of the number of the natural solutions (using CM) of one type of homogeneous Diophantine equation for values $k \geq 1 + log_2 s$, where $k$ the degree of the equation is and $s$ is the number of variables. It was also found the upper bound of the number of the natural solutions of the homogeneous algebraic Diophantine equation with integer coefficients with a small number of variables. Author investigated the relations of upper and lower estimates of the number of natural solutions of homogeneous Diophantine equation with integer coefficients diagonal form in the paper.

## 1. INTRODUCTION

The study of Diophantine equations are carried out in different directions: to find conditions for the existence of solutions of Diophantine equations, finding solutions of the Diophantine equation, to estimate the number of solutions of Diophantine equations. The paper focuses primarily on the last line of research.

Estimating of the number of solutions of algebraic Diophantine equations are considered over finite and infinite fields and rings. Estimating of the number of solutions of algebraic Diophantine equations over finite fields is carried out very successfully. It was found the estimate of the number of solutions of the homogeneous algebraic Diophantine equation diagonal form with coefficients from finite residues modulo [1].

The situation is worse with estimates of the number of solutions of algebraic Diophantine equations over infinite fields of rational numbers and the rings of integers, as these estimates are mainly qualitative in nature. Faltings [2] proved the Mordell conjecture about the finiteness of the number of rational points on the curve of genus g> 1.

Tue method (based on the results of Diophantine approximations) were developed in the works of Siegel [3], set the famous theorem on the finiteness of the number of integer points on a curve of genus $g \geq 1$.

Tue method has been extended by Schmidt [4] in the case of several variables, which allowed him to receive a multi-dimensional generalization of the Tue results about the finite number of integer solutions of Diophantine regulatory Tue equation. However, this result is also not effective.

Many methods of the theory of Diophantine equations, including Tue method, are not effective, thus do not allow a quantitative analysis of the solutions of Diophantine equations, therefore developed an effective method based on the use of lower estimates of logarithms of algebraic numbers.

Currently, quantitative estimates for the integer solutions of some classical Diophantine equations obtained by this method. Baker [5] made an effective assessment of the Tue equation: $f(x, y) = m$. A similar effective analysis for more general equation Tue-Mahler: $f(x, y) = mp_1^{x_1}...p_s^{x_s}$ (where $p_1,..., p_s$ are fixed numbers) was carried out in [6]. Another class of Diophantine equations (that allow effective analysis) constitutes superelliptic equation: $y^s = f(x)$, where $s \geq 2$, $f$ is an integer function of degree $n \geq 3$. Effective analysis of this equation was conducted by Baker [7]. This result was significantly strengthened [8]. Thus, the effective analysis is done only for algebraic Diophantine equations of two variables.

The circle method of Hardy and Littlewood (CM) [9] is independently from these studies. It allows making the top estimations of the number of natural solutions of the homogeneous algebraic Diophantine diagonal equation with integer coefficients. This method has been significantly strengthened by Vinogradov [10]. KM was used [11] to obtain top estimates of the number of natural solutions of algebraic Diophantine equations with integer coefficients. However, KM allows making these estimates, when the number of variables is much greater than the degree of homogeneous Diophantine equation.

Heath-Brown [14], [15] uses the methods of exponential sums and the determinant for estimating the number of integer (natural) solutions of homogeneous Diophantine diagonal equations with a small number of variables. However, these methods give a lower accuracy for considered classes of Diophantine equations in the paper.

In general, the lower estimate of the number of integers (positive) a solution of Diophantine equations is trivial (no solutions). We consider different cases of homogeneous Diophantine diagonal equations when this estimate is not trivial.

## 2. ESTIMATING THE NUMBER OF SOLUTIONS OF THE SIMPLIEST HOMOGENEOUS ALGEBRAIC DIOPHANTINE EQUATIONS

First, we consider lower estimations for the number of natural solutions of the simplest homogeneous Diophantine equations with $s$ variables. These estimates will be carried out in a hypercube with dimension $s$ and side $N$ - $A^s$, where $A = 1, ..., N$.

Let consider the simplest case:

$$x^k - y^k = 0.$$ (2.1)

It is clear that equation (2.1) has the solution $x = y$ and therefore it has $R_2^+(N) = N$ natural solutions (regardless of the value $k$) in $A^2$.

Homogeneous algebraic equation:

$$x_1^k - y_1^k + x_2^k - y_2^k = 0$$ (2.2)

has solutions $x_1 = y_1$, $x_2 = y_2$ and therefore it has $N^2$ natural solutions of equation (2.2) in $A^4$.

However, in contrast to (2.1), equation (2.2) has additional natural solutions, if $x_1 \neq y_1, x_2 \neq y_2$. For example, if $x_1 = a, y_1 = b$ then the symmetric natural solutions will be $x_2 = b, y_2 = a$, where $a, b$ - natural numbers. If $a > b$, then there is $C_N^2$ such natural solutions in $A^4$. If $b > a$, then it will be the same number of natural solutions in $A^4$. Thus, the equation (2.2) will have $2C_N^2 = N(N-1)$ additional salutations Therefore, the total number of natural solutions of the equation (2.2) (in $A^4$) will be:

$$R_4^+(N) = N^2 + N(N-1) = 2N^2 - N = O(N^2)$$ (2.3)

Now we consider the homogeneous equation in the form:

$$x_1^k - y_1^k + x_2^k - y_2^k + x_3^k - y_3^k + x_4^2 - y_4^k = 0. \qquad (2.4)$$

Equation (2.4) has the solutions: $x_1 = y_1, x_2 = y_2, x_3 = y_3, x_4 = y_4$ and therefore it has $N^4$ such natural solutions. However, the equation (2.4) has other additional solutions. This equation has (as was shown for equation (2.2)), $N(N-1)$ such other natural solutions, if $x_1 = y_1 = a, x_2 = y_2 = b$ ($a, b$ are fixed positive integers) and $x_3 \neq y_3, x_4 \neq y_4$. Similar number of natural solutions would be in case $x_1 = y_1 = a, x_3 = y_3 = b$, where $a, b$ are fixed natural numbers and $x_2 \neq y_2, x_4 \neq y_4$. There will be $C_4^2$ such cases. Therefore, it will be $C_4^2 N(N-1)$ additional natural solutions.

Now it is necessary to take into account that each natural numbers $a, b$ individually can have $N$ natural values. Therefore, the number of additional natural solutions will be $C_4^2 N(N-1)N^2$. We must also take into account additional natural solutions of the equation (2.4) in the case: $x_1 \neq y_1, x_2 \neq y_2, x_3 \neq y_3, x_4 \neq y_4$. Such additional natural solutions (as shown for equation (2.2)) will be: $(2C_N^2)(2C_N^2) = N^2(N-1)^2$. Thus, the equation (2.4) will have: $C_4^2 N(N-1)N^2 + N^2(N-1)^2$ additional natural solutions.

Therefore, homogeneous equation (2.4) will have the following number of natural solutions in $A^8$:

$$R_8^+(N) = N^4 + C_4^2 N(N-1)N^2 + N^2(N-1)^2 = 8N^4 - 8N^3 + N^2 = O(N^4). \qquad (2.5)$$

Now we consider the general case of an homogeneous equation of the form:

$$x_1^k - y_1^k + x_2^k - y_2^k + ... + x_{2^{j-1}}^k - y_{2^{j-1}}^k = 0 \qquad (2.6)$$

Equation (2.6) has the solution: $x_1 = y_1, x_2 = y_2, ..., x_{2^{j-1}} = y_{2^{j-1}}$ and therefore it has $N^{2^{j-1}}$ such natural solutions in $A^{2^j}$. The equation has also additional natural solutions in $A^{2^j}$. It has $N(N-1)$ additional natural solutions, as the equation (2.5), if $x_1 = y_1 = a_1, ..., x_{2^{j-1}-2} = y_{2^{j-1}-2} = a_{2^{j-1}-2}$ and $x_{2^{j-1}-1} \neq y_{2^{j-1}-1}, x_{2^{j-1}} \neq y_{2^{j-1}}$, where $a_1, ..., a_{2^{j-1}-2}$ are fixed natural numbers. It will be $C_{2^{j-1}}^2$ such cases, when 2 variables are not equal, and the remaining pairs are equal. So it will be $C_{2^{j-1}}^2 N(N-1)$ such additional natural solutions. We must take into account, that the natural numbers $a_1, ..., a_{2^{j-1}-2}$ (each separately) can have $N$ natural values. Therefore, the number of such additional natural solutions will be $C_{2^{j-1}}^2 N(N-1)N^{2^{j-1}-2}$.

Now it is necessary to consider the case when the 4 variables are pairwise equal, and the remaining pairs are equal. We obtain $C_{2^{j-1}}^4 N^2(N-1)^2 N^{2^{j-1}-4}$ additional natural solutions in this case etc. We must also consider the case when all variables are pairwise not equal. We obtain: $N^{2^{j-2}}(N-1)^{2^{j-2}}$ additional natural solutions in this case.

Consequently, the homogeneous equation (2.6) will have the following general number of natural solutions:

$$R_{2^j}^+(N) = N^{2^{j-1}} + \Sigma_{i=1}^{j-1} C_{2^{j-1}}^{2^i} N^{2^{i-1}} (N-1)^{2^{i-1}} N^{2^{j-1}-2^i} = N^{2^{j-1}} + \Sigma_{i=1}^{j-1} C_{2^{j-1}}^{2^i} (N-1)^{2^{i-1}} N^{2^{j-1}-2^{i-1}} = O(N^{2^{j-1}}) . (2.7)$$

Formulas (2.3), (2.5), (2.7) do not account all solutions, but only symmetric, so that these formulas do not depend from value $k$. However, the solutions may not be only being symmetric. For example, the equation (2.2) at the value $k = 4$ has solutions: $x_1 = 133, x_2 = 59, x_3 = 134, x_4 = 158$. Therefore, these estimates are lower bound for the number of natural solutions homogeneous Diophantine diagonal equation.

Based on (2.7) we obtain the following lower bound for the number of natural solutions of the equation (taking into account that the value is equal to $s = 2^j$ in equation (2.4)):

$$R_s^+(N) \geq O(N^{s/2}) . \tag{2.8}$$

It will be shown that (2.8) is the upper limit of the lower estimates for the number of natural solutions of the homogeneous algebraic diagonal form equation in $A^s$.

An upper bound of the number of natural solutions of the simplest homogeneous Diophantine diagonal equations will do with the help of the CM.

Let consider the homogeneous equation:

$$x_1^k - y_1^k + x_2^k - y_2^k + ... + x_{2^{j-1}}^k - y_{2^{j-1}}^k = 0, \tag{2.9}$$

where values $1 \leq x_i, y_i \leq N$.

It is well known Hua's lemma [9]. Suppose $1 \leq j \leq k$, then $\int_0^1 | f(x) |^{2^j} dx << N^{2^j - j + \epsilon}$, where $f(x) = \Sigma_{m=1}^N e^{2\pi x m^k}$.

Number of natural solutions of the equation (2.9) is defined by the integral [9]:

$$R_{2^j}^+(N) = \int_0^1 | f(x) |^{2^j} dx. \tag{2.10}$$

Therefore, based on Hua's lemma we have the following upper bound for the number of natural solutions of the equation (2.9):

$$R_{2^j}^+(N) << N^{2^j - j + \epsilon}, \tag{2.11}$$

where value $\epsilon$ is a small positive real number and $1 \leq j \leq k$. The case $j > k$ will consider when solving the following equation.

Let us consider a more general homogeneous equation:

$$x_1^k - y_1^k + x_2^k - y_2^k + ... + x_s^k - y_s^k = 0, \tag{2.12}$$

where values $1 \leq x_i, y_i \leq N$.

Based on [9] the number of natural solutions of the equation (2.12) defined by the integral:

$$R_{2s}^+(N) = \int_0^1 | f(x) |^{2s} dx. \tag{2.13}$$

Thus it is necessary to make an estimate of the integral (2.13) for the estimation of the number of natural solutions of the equation (2.12).

The number of representations of the number $n$ as a sum of $s$ members of $k$ degrees of natural numbers (based on [9]) is defined by the formula:

$$R_s^p(n) = \int_0^1 f(x)^s e^{-2\pi x n} dx = \Gamma(1 + 1/k)\Gamma(s/k)^{-1} n^{s/k - 1} \Sigma(n) + O(n^{s/k - 1 - \delta}), \tag{2.14}$$

where $f(x) = \Sigma_{m=1}^{[n^{1/k}]} e^{2\pi x m^k}$, $\Gamma$-gamma function, $\Sigma(n)$- the sum of the series, $\delta$ - a small positive real number and $s > 2^k$.

Based on (2.14) let us make an asymptotic upper bound of value $R_s^p(n)$ (2.14):

$$| R_s^p(n) | = \int_0^1 | f(x) |^s | e^{-2\pi xn} | dx = | \Gamma(1+1/k)\Gamma(s/k)^{-1} n^{s/k-1}\Sigma(n) + O(n^{s/k-1-\delta}) | << n^{s/k-1+\epsilon}. \quad (2.15)$$

Based on (2.15): $| e^{-2\pi xn} | = 1$, $N = [n^{1/k}]$ and $\Gamma(1+1/k)\Gamma(s/k)^{-1} < 1, \Sigma(n) << n^\epsilon$ we obtain the following estimate:

$$\int_0^1 | f(x) |^s dx << n^{s/k-1+\epsilon} = N^{s-k+\epsilon}, \quad (2.16)$$

for values $s > 2^k$.

It follows the required estimate (based on (2.16)):

$$\int_0^1 | f(x) |^{2s} dx << N^{2s-k+\epsilon}, \quad (2.17)$$

for values $s > 2^{k-1}$.

Formula (2.17) gives the upper bound for the number of natural solutions of the equation (2.12) for values $s > 2^{k-1}$. Therefore (according to the formula (2.17)) we can estimate the number of natural solutions of the equation (2.12) for values $k < 1 + log_2 s$ for a fixed value $s$.

It is interesting to get an estimate of the number of natural solutions of the equation (2.12) for values $k \geq 1 + log_2 s$. We use the fact that every even number can be represented as the sum of two degrees with natural indicators $j_1 > ... > j_t \geq 1$:

$$2s = 2^{j_1} + ... + 2^{j_t}, \quad (2.18)$$

Based on (2.18) and Cauchy-Schwarz inequality, we obtain:

$$\int_0^1 | f(x) |^{2s} dx \leq (\int_0^1 | f(x) |^{2^{j_1+1}} dx^{1/2}) \cdot (\int_0^1 | f(x) |^{2^{j_2+2}} dx)^{1/4} \cdot \cdot (\int_0^1 | f(x) |^{2^{j_t+t}} dx)^{1/2^{t-1}}. \quad (2.19)$$

Based on (2.19) and Hua's lemma we obtain:

$$(\int_0^1 |f(x)|^{2^{j_t+i}} dx)^{1/2^i} << N^{2^{j_t} - (j_t+i)/2^i + \epsilon/2^i}. \quad (2.20)$$

Thus, based on (2.19), (2.20) we obtain (for values $k \geq j_1 = 1 + log_S$):

$$R_{2s}^+(N) = \int_0^1 | f(x) |^{2s} dx << N^{\Sigma_{i=1}^{t-1}(2^{j_i} - (j_i+i)/2^i)/2^i + \epsilon/2^i + 2^{j_t+1} - (j_t+t)/2^{t-1} + \epsilon/2^{t-1}} = N^{2s+\epsilon - \Sigma_{i=1}^{t-1}(j_i+i)/2^i + (j_t+t)/2^{t-1}}. \quad (2.21)$$

For example, based on (2.21), we obtain (for values $s = 15, 2s = 30 = 2^4 + 2^3 + 2^2 + 2$ ($j_1 = 4, j_2 = 3, j_3 = 3, j_4 = 1$) and values $k \geq j_1 = 4$):

$$R_{30}^+(N) << N^{30+\epsilon-5} = N^{25+\epsilon} . \tag{2.22}$$

Based on (2.17) let us determine the number of solutions of the equation (for values $k < 4$):

$$\int_0^1 |f(x)|^{30} \, dx << N^{30-k+\epsilon} . \tag{2.23}$$

It is easy to see (comparing (2.23) with (2.22)), that the number of solutions of the equation for smaller values $k$ is more. That is true.

## 3. ESTIMATIONS OF THE NUMBER OF NATURAL SOLUTIONS OF HOMOGENEOUS DIOPHANTINE DIOGANAL EQUATION WITH LARGE NUMBER OF VARIABLES

Now we consider the homogeneous diagonal equation with integer coefficients:

$$a_1 x_1^k + ... + a_s x_s^k = 0 , \tag{3.1}$$

where $k$ is a positive integer, and all values $a_i$ are integers.

Fermat's equation (by the way) belongs to the class of equations:

$$x_1^k + x_2^k - x_3^k = 0 . \tag{3.2}$$

It was proved the following theorem using CM in [9].

Пусть $k \geq 2$ и $s_0$, как в теореме Виноградова (см. ниже) и пусть $s \geq min(s_0, 2^k + 1)$ и $s \geq 4k^2 - 4k + 1$. Предположим также, что когда $k$ четное, то не все целые числа $a_1, ... a_s$ имеют один знак. Тогда уравнение (3.1) имеет нетривиальные решения в целых числах $x_1, ... x_s$ ..

Suppose that $k \geq 2$ and $s_0$ as in Vinogradov's theorem (see below) and let $s \geq min(s_0, 2^k + 1)$ and $s \geq 4k^2 - 4k + 1$. We also assume when $k$ an even is, then not all integers $a_1, ... a_s$ have the same sign. Then the equation (3.1) has non-trivial solutions in integers $x_1, ... x_s$ .

Note. One can also be considered that not all $a_i$ have the same sign for odd values $k$. If it is necessary, one can replace the unknown $x_i$ on $-x_i$.

Vinogradov's theorem gives values $s_0 < 2^k + 1$ for values $k \geq 11$ and values $s_0 \geq 2^k + 1$ for values $k < 11$. Therefore, it must be satisfied: $s \geq 2^k + 1$ and $s \geq 4k^2 - 4k + 1$ for values $k < 11$. For example, it is performed $s \geq 15$ for value $k = 2$ and $s \geq 34$ for value $k = 3$. It is performed $2^k + 1 > 4k^2 - k + 1$ only for values $k \geq 8$, thus enough $s \geq 2^k + 1$.

Based on CM we find an upper estimate for the number of natural solutions of the equation (3.1) in a hypercube with side $N$:

$$R_s^+(N) = \Sigma \cdot J(N) + O(N^{s-k+\varepsilon}), \tag{3.3}$$

where $\Sigma$ is the sum of a convergent series, and

$$N^{s-k} << J(N) << N^{s-k+\varepsilon}. \tag{3.4}$$

Based on (3.3) and (3.4) we obtain:

$$R_s^+(N) << N^{s-k+\varepsilon}. \tag{3.5}$$

The above theorem does not mean that the equation with a lower value can't have an infinite number of solutions. For example, let consider the equation:

$$x_1^k + x_2^k - Dx_3^k = 0 \tag{3.6}$$

Equation (3.6) has an infinite number of solutions, at a value $D = 2n^k$ which is equal to: $x_1 = x_2 = nx_3$.

Other examples of homogeneous diagonal equations with an infinite number of solutions for small values $s$ are homogeneous equations previously discussed (2.1) and (2.2).

## 4. ESTIMATION OF THE NUMBER OF NATURAL SOLUTIONS OF HOMOGENEOUS DIAGONAL DIOPHANTINE EQUATION WITH A SMALL NUMBER OF VARIABLES

KM applies to the determination of the number of natural solutions of homogeneous diagonal Diophantine equation at large values $s$ in comparison with $k$ (more precisely see the theorem above). Thus it is interesting to determine the number of natural solutions of the equation (3.1) for small values $s$ using other methods.

Let consider natural solutions of homogeneous algebraic diagonal Diophantine equation:

$$a_{x1}^{k} + ... + a_{s}x_{s}^{k} = 0,$$ (4.1)

where $k$ is natural number and factors $a_i$ are relatively prime integers. While at least one of the factors $a_i$ has a sign different from the other (this is not required if $k$ is odd). Any homogeneous diagonal Diophantine equation reduces to the equation (4.1), if it is divided by the greatest common divisor of its coefficients.

It is well known that if the equation (4.1) has at least one natural solution $x_{10}, ..., x_{s0}$, then it has an infinite number of natural solutions. We assume performing the following condition: $x_{10} \leq ... \leq x_{s0}$. It is performed following lower bound for the number of natural solutions of the equation (4.1) in the hypercube with side - $N$ :

$$R_s(N) \geq N / inf(x_{s0}),$$ (4.2)

if the equation has at least one natural solution.

Let denote $D(f)$ as algebraic variety of the Diophantine equation $f(x_1, ..., x_s) = 0$. We consider $f$ as a polynomial with integer coefficients. We also consider solutions of the Diophantine equation $f(x_1, ..., x_s) = 0$ in $(Z^+)^s, s \geq 1$ (a direct $s$ product of the set of natural numbers).

We prove a simple assertion, which is applicable not only to the homogeneous algebraic Diophantine equations, but to a wider class of Diophantine equations.

Assertion 1

Let the polynomial $f(x_1, ...x_r, ...x_s)$ satisfies the condition:

$$f(x_1, ...x_r, ...x_s) = f_1(x_1, ...x_r) + f_2(x_{r+1}, ...x_s)$$ (4.3)

in $(Z^+)^s$ then:

1. If the value $f_1 \cdot f_2 > 0$ in the rest of area $(Z^+)^s$ then:

$$R_s(N) = R_r(N) \bullet R_{s-r}(N)$$ (4.4)

in the hypercube with the dimension - $s$ and side - $N$, where $R_i(N)$ is the number of solutions of the equation $f(x_1,..,x_i)=0$ in the area $(A)^s$.

2. If the value $f_1 \cdot f_2 < 0$ in the rest of area $(Z^+)^s$ then:

$$R_s(N) \geq R_r(N) \cdot R_{s-r}(N), \qquad (4.5)$$

in the hypercube with the dimension - $s$ and side - $N$.

We prove the first 1. Let $D(f)$ is an algebraic variety of the equation $f(x_1,...x_r,...x_s)=0$, $D(f_1)$ is an algebraic variety of the equation $f_1(x_1,...x_r)=0$ and $D(f_2)$ is an algebraic variety of the equation $f_2(x_{r+1},...x_s)=0$. Then it is satisfied $f(x_1...x_r,...x_s)=f_1(x_1,...x_r)+f_2(x_{r+1},...x_s)=0$ for an algebraic variety $D(f)=D(f_1)\times D(f_2)$ ($\times$-direct product). It is performed $f_1 \cdot f_2 > 0$ in the rest of the area $(Z^+)^s$ ($f_1>0, f_2>0$ or $f_1<0, f_2<0$). Thus we do not have other solutions.

Now we prove 2. Let $D(f)$ is an algebraic variety of the equation $f(x_1,...x_r,...x_s)=0$, $D(f_1)$ is an algebraic variety of the equation $f_1(x_1,...x_r)=0$ and $D(f_2)$ is an algebraic variety of the equation $f_2(x_{r+1},...x_s)=0$. Then it is satisfied $f(x_1...x_r,..x_s)=f_1(x_1,...x_r)+f_2(x_{r+1},...x_s)=0$ for an algebraic variety $D(f)=D(f_1)\times D(f_2)$ ($\times$-direct product). It is performed $f_1 \cdot f_2 < 0$ in the rest of the area $(Z^+)^s$ ($f_1>0, f_2<0$ or $f_1<0, f_2>0$). Thus we can have other solutions of the equation $f=f_1+f_2=0$ and $R_s(N) \geq R_r(N) \cdot R_{s-r}(N)$.

The assertion 1 is easily offended by region integers $(Z)^s$ and inhomogeneous equations.

Here is an example using the assertion 1. Let us consider the inhomogeneous equation $f = x_1^2 + x_1 - x_2^3 + x_2 = 0$. This equation satisfies (4.3) for the functions $f_1 = x_1^2 + x_1 = 0$ and $f_2 = -x_2^3 + x_2 = 0$. The equation has two solutions in $(Z)$: $x_{11}=0, x_{12}=-1$. The equation $f_2 = -x_2^3 + x_2 = 0$ has three solutions in $(Z)$: $x_{21}=-1, x_{22}=0, x_{23}=1$. As we can see there are solutions: $(-1,-1),(-1,0),(-1,1),(0,-1),(0,0),(0,1)$ on the algebraic variety $D(f_1) \times D(f_2)$ in $(Z)^2$. It is performed $f_1 \cdot f_2 < 0$ in the rest of the area $(Z)^2$, so it is satisfied the second part of the assertion 1. Based on (4.5) the equation $x_1^2 + x_1 - x_2^3 + x_2 = 0$ can have more solutions than a direct product of the algebraic variety of the equations: $f_1 = x_1^2 + x_1 = 0$ and $f_2 = -x_2^3 + x_2 = 0$. We have 6 solutions on the algebraic variety $D(f_1) \times D(f_2)$ in $(Z)^2$. However, it is more solutions of the equation $f = x_1^2 + x_1 - x_2^3 + x_2 = 0$ in $(Z)^2$. It is easy to see that the solution $(2,2)$ also satisfies the equation $f = x_1^2 + x_1 - x_2^3 + x_2 = 0$.

Let us look at another example. In [12] asks the question about the number of natural solutions of non-homogeneous equation $x^2 + y^2 = z^2 + 1$. Based on identity $(n^2 + n - 1)^2 + (2n + 1)^2 = (n^2 + n + 1)^2 + 1$ the equation has an infinite number of solutions.

Let consider more general equation: $x^l + y^m = z^m + a^l$, where $l, m, a$ are arbitrary natural numbers. It is not possible to answer the question - whether the inhomogeneous equation has an infinite number of solutions based the identity. However, this can be done based on assertion 1. We write this equation in the form: $x^l - a^l + y^m - z^m = 0$. Let us denote: $f_1 = x^l - a^l$, $f_2 = y^m - z^m$. The equation $f_1 = 0$ has one natural solution $x = a$ and the equation $f_2 = 0$ has an infinite number of natural solutions $y = z$. This equation has $R_2(N) = N$ natural solutions in $(A)^2$. It is performed $f_1 f_2 < 0$, as $f_1, f_2$ have different signs in the rest area. Thus based on assertion 1 an algebraic variety $D(f_1) \times D(f_2)$ is infinitely and the lower estimate of the number of natural solutions of the equation $f = 0$ will be $R_3^+(N) \geq N$ in $(A)^3$.

Consequence 1. The algebraic Diophantine equations in two variables $f(x_1) + f(x_2) = 0$ has a finite number of integer salutations does not exceed the product of powers of polynomials - $\deg f_1 \cdot \deg f_2$, if the conditions $f_1 \cdot f_2 > 0$ are fulfilled in the rest area $(Z)^2$.

This consequence is easy to communicate in case of $s$ variables.

Consequence 2

Let a polynomial $f(x_1,...,x_s)$ satisfies (4.3), i.e. $f(x_1,...,x_s) = f_1(x_1,...,x_m) + f_2(x_{m+1},...,x_s)$ and the homogeneous equation $f_1(x_1,...x_m) = 0$ has $R_m^+(N)$ natural solutions and homogeneous equation $f_2(x_{m+1},...x_s) = 0$ has $R_{s-m}^+(N)$ natural solutions. Then for the number of natural solutions of the homogeneous equation $c_1 f_1 + c_2 f_2 = 0$ (where $c_1, c_2$ are integers) we have the estimate $R_s^+(N) \geq R_m^+(N) \cdot R_{s-m}^+(N)$.

Definition. Polynomial $a_1 x_1^k + \ldots + a_s x_s^k$ (where $a_i$ - integers $(1 \le i \le s)$) is called partially symmetrical about variables $x_l, x_m$ if performed $a_l = a_m$.

It is necessary to take into account the permutations of solutions for determining the number of solutions of the homogeneous equation $a_1 x_1^k + \ldots + a_s x_s^k = 0$, where $a_i$ - integers $(1 \le i \le s)$ and $a_l = a_m$. Thus, if this equation has a natural solution $(x_{10}, \ldots x_{l0}, \ldots x_{m0}, \ldots x_{s0})$, it also has a solution $(x_{10}, \ldots x_{m0}, \ldots x_{l0}, \ldots x_{s0})$.

Now let us return to the consideration of the homogeneous Diophantine equation (4.1).

First, let us consider the simplest case $s = 2$, which we will need in the future. The necessary condition for the existence of solutions of the homogeneous Diophantine diagonal equation (as already mentioned) is different signs of the coefficients. Thus, let us consider the equation:

$$a_1 x_1^k - a_2 x_2^k = 0, \tag{4.6}$$

where $a_1, a_2, k$ are natural numbers.

Suppose, for definiteness $0 < a_1 \le a_2$, then the equation (4.6) is equivalent to the equation: $x_1^k - a_2 / a_1 x_2^k = (x_1 - (a_2 / a_1)^{1/k} x_2) \cdot P(x_1, x_2)$, where the polynomial $P(x_1, x_2)$ has no natural zero. Therefore, the equation (35) has the natural solution in the straight line:

$$x_1 = (a_2 / a_1)^{1/k} x_2, \tag{4.7}$$

if the value $(a_2 / a_1)^{1/k}$ is rationally.

If $x_2 = 1$ then the value is $x_1 = (a_1 / a_2)^{1/k} \ge 1$. Therefore, based on (4.7) the number of the natural solutions of the equation (4.6) is:

$$R_2(N) = N / ((a_2 / a_1)^{1/k}) \le N \tag{4.8}$$

in $(A)^2$.

Now we consider the homogeneous Diophantine equation (4.1) for the case $s = 3$:

$$a_1 x_1^k + a_2 x_2^k - a_3 x_3^k = 0, \tag{4.9}$$

where all $a_i$ are natural. The equation: $a_1^{\ k} x_1^k - a_2 x_2^k - a_3 x_3^k = -(a_2 x_2^k + a_3 x_3^k - a_1 x_1^k) = 0$ reduces to the form (4.9).

Let us consider the equation (4.9) for the value $k = 1$:

$$a_1 x_1 + a_2 x_2 - a_3 x_3 = 0. \qquad (4.10)$$

Note that the equation $a_1 x_1 + a_2 x_2 + ... + a_{s-1} x_{s-1} - a_s x_s = 0$ ($s \geq 2$ and $a_i$ are natural numbers) has an infinite number of natural solutions in the general case.

Equation (4.10) corresponds to the plane passing through the origin at which the normal vector has a positive value on the axis of projection $x_1, x_2$ and a negative value of the projection on the axle $x_3$.

Suppose that $0 < a_1 \leq a_2 \leq a_3$, then the equation (4.10) can be written as:

$$x_3 = (a_1 / a_3) x_1 + (a_2 / a_3) x_2. \qquad (4.11)$$

Let determine the number of natural solutions of the equation (4.10). Suppose $x_1 = i, x_2 = t_i$, then based on (4.11) we obtain $x_3 = i(a_1 / a_2) + t_i(a_2 / a_3) \leq N$ and $t_i \leq (a_3 / a_2)(N - i(a_1 / a_2))$.

Therefore, the number of natural solutions of the equation (4.10) is equal to:

$$R_3(N) = \Sigma_i(t_i) = (a_3 / a_2)\Sigma_i(N - i(a_1 / a_2)) = (a_3 / a_2)(a_2 / a_1 - a_1 / 2a_2)N(N-1) = O(N^2). \quad (4.12)$$

Now we consider the equation (4.9) with the value $k = 2$:

$$a_1 x_1^2 + a_2 x_2^2 - a_3 x_3^2 = 0. \qquad (4.13)$$

Equation (4.13) corresponds to the cosine of the second order, the apex of which is located at the origin.

Let us suppose $a_1 \geq a_2 \geq a_3 > 0$. In this case equation (4.13) is equivalent to the equation:

$$x_3^2 = (a_1 / a_3) x_1^2 + (a_2 / a_1) x_2^2. \qquad (4.14)$$

Let us make an affine transformation of coordinates:

$$y_1 = (a_1 / a_2)^{1/2} x_1, \ y_2 = (a_2 / a_3)^{1/2} x_2, \ y_3 = x_3. \qquad (4.15)$$

If $(a_1/a_2)^{1/2}$ and $(a_2/a_3)^{1/2}$ are positive integers, then the affine coordinate transformation (4.15) carries the natural solution of the equation (4.14) in a natural solution of the equation:

$$y_1^2 + y_2^2 = y_3^2.\qquad(4.16)$$

Equation (4.16) is the equation of Fermat with the value $k = 2$ that it has an infinite number of solutions:

$$y_1 = (m^2 - n^2)l,\ y_2 = 2mnl,\ y_3 = (m^2 + n^2)l,\qquad(4.17)$$

where $m, n, l$ are natural numbers, while $m > n$.

We see that the maximum value has the variable $y_3$. Therefore, the number of the natural solutions is:

$$y_3 = (m^2 + n^2)l \le N\qquad(4.18)$$

in $(A)^3$.

Thus, it is necessary to find the number of points with the natural coordinates, satisfying the condition (4.18).

If $l = 1$ than we obtain: $m^2 + n^2 \le N$ from (4.18). Since $m > n$ we obtain the number of points with the natural coordinates within the sector above the main diagonal with a radius equal. $N^{1/2}$ Based on [13], we have such points:

$$\pi N/8 + O(N^{1/2}).\qquad(4.19)$$

Based on (4.19) (when $l = l_{max}$) we get inequality: $m^2 + n^2 \le N/l_{max}$. This is the number of points with natural coordinates within a sector above the main diagonal with a radius equal $(N/l_{max})^{1/2}$. Therefore, these points will be:

$$\pi N/8l_{\max} + O((N/l_{\max})^{1/2}).\qquad(4.20)$$

Based on (4.19), (4.20) we have the following upper estimate for the number of natural solutions of the equation (4.16):

$$R_3^+(N) \le \pi N/8(1 + 1/2 + \ldots + 1/l_{\max}) + O((N)^{1/2}).\qquad(4.21)$$

The function $1/k$ is strictly decreasing function and, therefore, based on the Euler-Maclaurin formula we get:

$$\Sigma_{k=1}^{lmax} 1/k = \int_{k=1}^{lmax} dx/x + C + O(1/l_{max}) = ln(l_{max}) + C + O(1/l_{max}), \qquad (4.22)$$

where $C$ is constant.

Substituting (4.22) into (4.21) and obtain:

$$R_3^+(N) \leq \pi N/8[ln(l_{max}) + C + O(1/l_{max})] + O((N)^{1/2}). \qquad (4.23)$$

Based on (4.23) and $l_{max} = N/5$ we obtain an asymptotic upper estimate of the number of natural solutions of the equation (4.16):

$$R_3^+(N) \leq \pi N ln(N)/8 + O(N). \qquad (4.24)$$

When $a_1 \geq a_2 \geq a_3 > 0$ and $(a_1/a_3)^{1/2}, (a_2/a_3)^{1/2}$ are natural (based on (4.24)) we have the following asymptotic upper bound for the number of natural solutions of the equation (4.13):

$$R_3^+(N) \leq \pi N ln(N)/8(a_1/a_3)^{1/2} + O(N) << N^{1+\epsilon}. \qquad (4.25)$$

If $k > 2$ and the equation $a_1 x_1^k + a_2 x_2^k - a_3 x_3^k = 0$ has at least one natural solution, then we have the following lower bound for the number of natural solutions of the equation (4.9) in $(A)^3$:

$$R_3^+(N) \geq N/sup_{i \leq 3} x_{i0}. \qquad (4.26)$$

On other hand, Fermat's equation $x_1^k + x_2^k - x_3^k = 0$ has no natural solutions at $k > 2$.

Now we consider the homogeneous equation (4.1) with the value $s = 4$:

$$a_1 x_1^k - a_2 x_2^k + a_3 x_3^k - a_4 x_4^k = 0, \qquad (4.27)$$

where all $a_i$ are natural numbers.

Let us use the assertion proved above. Consider polynomials $f_1 = a_1 x_1^k - a_2 x_2^k$ and $f_2 = a_3 x_3^k - a_4 x_4^k$. We assume (without loss of generality), that $a_2 \geq a_1, a_4 \geq a_3$. Based on the (4.8) the equation $f_1 = a_1 x_1^k - a_2 x_2^k = 0$ has the following number of natural solution in $(A)^2$: $R_2(N) = N/(a_2/a_1)^{1/k}$, if $(a_2/a_1)^{1/k}$ is a rational number. Based on the (4.8) the equation $f_2 = a_3 x_3^k - a_4 x_4^k = 0$ has the following number of natural solution in $(A)^2$: $R_2(N) = N/(a_4/a_3)^{1/k}$, if $(a_4/a_3)^{1/k}$ is a rational number. It is performed the condition $f_1 \cdot f_2 < 0$ in the rest part of the area $(Z^+)^4$, so (based the second part of assertion 1) under the above conditions, the following lower bound of l number of natural solutions of the equation (4.27) is:

$$R_4^+(N) \geq N^2 / (a_2 a_4 / a_1 a_3)^{1/k}.\tag{4.28}$$

If $a_2 \geq a_1, a_4 \geq a_3$, then $a_2 a_4 \geq a_1 a_3$. Consequently, a maximum of (4.28) is reached at values $a_1 = a_2, a_3 = a_4$ and it is equal to:

$$R_4^+(N) \geq N^2.\tag{4.29}$$

On the other hand, the maximum upper bound of the estimation of the equation (4.27) is achieved at values $a_1 = a_2 = a_3 = a_4$ and corresponds to the equation:

$$x_1^k - x_2^k + x_3^k - x_4^k = 0.\tag{4.30}$$

Based on Hua's lemma the upper estimate of the number of natural solutions of the equation (4.30) is equal to:

$$R_4^+(N) << N^{2+\epsilon}.\tag{4.32}$$

where $\epsilon$ is a small real number (cf. (4.29)).

Now we consider the homogeneous algebraic Diophantine equation with an arbitrary odd value $s$:

$$a_1 x_1^k - a_2 x_2^k + \ldots + a_{s-1} x_{s-1}^k - a_s x_s^k = 0.\tag{4.33}$$

Assertion 2

Suppose $a_2 \geq a_1, \ldots, a_s \geq a_{s-1}$ in the equation (4.33), then the lower estimate of the number of natural solutions of the equation (4.33) (if $(a_2 a_4 \ldots a_s / a_1 a_3 \ldots a_{s-1})^{1/k}$) is a positive integer) is equal to:

$$R_s^+(N) \geq N^{s/2} / (a_2 a_4 \ldots a_s / a_1 a_3 \ldots a_{s-1})^{1/k}.\tag{4.34}$$

Proof. We perform the proof by induction on values $s$.

Based on (4.28) holds the formula (4.34) for the value $s = 4$.

Suppose that (4.34) holds for the value $s = 2t$, and check that it is carried out for the value $s = 2t + 2$.

Take polynomials $f_1 = a_1 x_1^k - a_2 x_2^k + .... + a_{2t-1} x_{2t-1}^k - a_{2t} x_{2t}^k$ and $f_2 = a_{2t+1} x_{2t+1}^k - a_{2t+2} x_{2t+2}^k$.

Based on the second part of Assertion1 it is performed:

$$R_{2t+2}^+(N) \geq N^{2t/2} / (a_2 a_4 ... a_{2t} / a_1 a_3 ... a_{2t-1})^{1/k} \cdot N^{2/2} / (a_{2t+2} / a_{2t+1})^{1/k} = N^{2t+2/2} / (a_2 a_4 ... a_{2t+2} / a_1 a_3 ... a_{2t+1})^{1/k},$$

as the condition $f_1 \cdot f_2 < 0$ performed in the rest part of the field $(Z^+)^{2t+2}$.

Since the condition: $a_2 \geq a_1, ..., a_s \geq a_{s-1}$, then $a_2 a_4 ... a_s / a_1 a_3 ... a_{s-1} \geq 1$. Therefore, the maximum lower bound (4.34) is achieved (at the value $a_2 a_4 ... a_s / a_1 a_3 ... a_{s-1} = 1$):

$$R_s^+(N) \geq N^{s/2}. \tag{4.35}$$

Maximum upper estimate of the number of natural solutions of the equation (4.31) is achieved at values $a_1 = ... = a_s$ that correspond to the equation:

$$x_1^k - x_2^k + ... + x_{s-1}^k - x_s^k = 0. \tag{4.36}$$

Based on Hua's lemma the upper estimate of the number of natural solutions of the equation (4.36) is equal to:

$$R_s^+(N) << N^{s/2+\epsilon}. \tag{4.37}$$

where $\epsilon$ is a small real number (cf. (4.35)).

Let us consider the homogeneous equation (4.1) with the value $s = 5$:

$$a_1 x_1^k + a_2 x_2^k + a_3 x_3^k + a_4 x_4^k - a_5 x_5^k = 0. \tag{4.38}$$

If the equation (4.38) has at least one natural root $(x_{10}, x_{20}, x_{30}, x_{40}, x_{50})$, then (regardless of the value $k$) the lower estimates of the number of natural solutions of the equation (4.38) is equal to:

$$R_5^+(N) \geq N / sup_{i \leq 5} x_{i0}. \tag{4.39}$$

For example, the equation $x_1^5 + x_2^5 + x_3^5 + x_4^5 - x_5^5 = 0$ has solutions $x_1 = 27, x_2 = 84, x_3 = 110, x_4 = 133, x_5 = 144$.

Therefore, based on (4.39) (taking into account the partial symmetry of the polynomial in the variables $x_1, x_2, x_3, x_4$ and the values of inequality $x_1 = 27, x_2 = 84, x_3 = 110, x_4 = 133$) the following lower bound for the number of natural solutions of the equations is true: $R_5^+(N) \geq 4! N / 144$.

Now we consider another algebraic equation with $s = 5$:

$$a_1 x_1^k + a_2 x_2^k - a_3 x_3^k + a_4 x_4^k - a_5 x_5^k = 0, \qquad (4.40)$$

where all $a_i$ - natural numbers.

First, we consider the equation (4.40) in the case $k = 1$. Let $f_1 = a_1 x_1 + a_2 x_2 - a_3 x_3 = 0$. Based on (4.12) the equation (under the assumption $0 < a_1 \leq a_2 \leq a_3$) has the number of natural solutions: $R_3^+(N) = a_3 / a_2 (a_2 / a_1 - a_1 / 2a_2) N(N-1)$ in $(A)^3$ and the equation $f_2 = a_4 x_4 - a_5 x_5 = 0$ (under the assumption $0 < a_4 \leq a_5$) has the number of natural solutions: $R_2^+(N) = N / (a_5 / a_4)$.

Based on accretion 1 (for values $f_1 f_2 < 0$ in $(Z^+)^5$), the equation (4.40) (under the assumptions above) has the number of natural solutions at least:

$$R_5^+(N) \geq a_3 a_4 / a_1 a_5 (a_2 / a_1 - a_1 / 2a_2) N^2 (N-1) = O(N^3). \qquad (4.41)$$

Now consider the equation (4.40) in the case $k = 2$. If $a_1 \geq a_2 \geq a_3$ and $(a_1 / a_3)^{1/2}$ is natural number, the equation $f_1 = a_1 x_1^2 + a_2 x_2^2 - a_3 x_3^2 = 0$ (based on (4.25)) has the following number of natural solutions: $R_3^+(N) = \pi N (a_1 / a_3)^{1/2} \ln(N) / 8$.

Based on (4.8) (if $(a_5 / a_4)^{1/2}$ is natural number) the equation $f_2 = a_4 x_4^2 - a_5 x_5^2 = 0$ has the following number of natural solutions in $(A)^2$: $R_2(N) = N / (a_5 / a_4)^{1/2}$.

Equation (4.40) (at the value $k = 2$) can be written as $f = f_1 + f_2 = 0$ and apply assertion 1 for the case $f_1 f_2 > 0$ (under the assumptions mentioned above), we obtain the following lower bound for the number of its natural solutions:

$$R_5^+(N) \geq \pi N^2 \ln(N) / 8(a_1 a_5 / a_3 a_4)^{1/2} = O(N^2 \ln(N)). \qquad (4.42)$$

Let us consider the equation (4.40) in the case $k > 2$. If the equation $f_1 = a_1 x_1^k + a_2 x_2^k - a_3 x_3^k = 0$ has at least one solution $x_{10}, x_{20}, x_{30}$, then it is performed the following lower bound of the estimation for the number of natural solutions of the equation: $R_3^+(N) \geq N / sup_{i \leq 3}(x_{i0})$. Based on (4.8) (if $(a_5 / a_4)^{1/k}$ is natural number) the equation $f_2 = a_4 x_4^k - a_5 x_5^k = 0$ has the following number of natural solutions in $(A)^2$: $R_2(N) = N / (a_5 / a_4)^{1/k}$.

The homogeneous equation (4.40) can be written in the form $f = f_1 + f_2 = 0$, so based on assertion 1 (and in the case $f_1 f_2 > 0$ and the above-mentioned assumptions for values $k > 2$) we obtain the following lower bound for the number of natural solutions of the equation:

$$R_5^+(N) \geq N^2 / (a_4 / a_5)^{1/k} sup_{i \leq 3}(x_{i0}) = O(N^2). \qquad (4.43)$$

## 5. LOWER ESTIMATE OF THE NUMBER OF NATURAL SOLUTIONS OF HOMOGENEOUS DIOPHANTINE DIAGONAL EQUIOPHAHTINE WITH ARBITRARY NUMBER OF VARIABLES

Let us consider the class of homogeneous irreducible algebraic Diophantine diagonal equations with an arbitrary number of variables in $(A)^s$:

$$c_1 x_1^k + ... + c_s x_s^k = 0, \qquad (5.1)$$

where $c_i (1 \leq i \leq s)$ are integers, and $k$ is a natural number.

Equation (5.1) can be written as:

$$f = f_1 + f_2 = 0 \qquad (5.2)$$

where $f_1 = a_1 x_1^k + .. + a_{m-1} x_{m-1}^k - a_m x_m^k = 0$ and $f_2 = a_{m+1} x_{m+1}^k - a_{m+2} x_{m+2}^k + ... + a_{s-1} x_{s-1}^k - a_s x_s^k = 0$.

Functions $f_1$ and $f_2$ may be omitted in a particular case.

Note that the number of monomials with negative coefficients in the equation (5.2) does not exceed $s / 2$. Otherwise, we will consider the equivalent equation $-f = 0$.

We group the monomials of the equation (5.2) so that based on assertion 1 to obtain a more effective lower bound of the number of natural solutions. We must select functions $f_1$, $f_2$ in such way (if it possible), that equations $f_1 = 0$ and $f_2 = 0$ have solutions in positive integers.

For example, we set in the equation:

$$x_1^2 + x_2^2 - 4x_3^2 + x_4^2 - 2x_5^2 = 0 \tag{5.3}$$

$f_1 = x_1^2 + x_2^2 - 2x_5^2 = 0$ and $f_2 = x_4^2 - 4x_3^2 = 0$.

Equation $f_1 = 0$ has the solution $x_1 = x_2 = x_5$ and $N$ natural solutions in $(A)^3$. Equation $f_2 = 0$ has the solution $x_4 = 2x_3$ in $(A)^2$ and $N/2$ natural solutions. Based on assertion 1 we obtain the following effective lower bound for the number of natural solutions of the equation (5.3) in the case $f_1 f_2 > 0$: $R_5^+(N) \geq N^2/2$.

If we take $f_1 = x_1^2 + x_2^2 - 4x_3^2 = 0$, then this equation has $N/2$ natural solutions in $(A)^3$. However, the equation $f_2 = x_4^2 - 2x_5^2 = 0$ does not have natural solutions in $(A)^2$ in this case, so based on assertion 1 (for values $f_1 f_2 < 0$) we obtain inefficient lower estimate for the number of natural solutions of the equation (5.3): $R_5^+(N) \geq 0$.

Now we consider separately the equation:

$$f_1 = a_1 x_1^k + .. + a_{m-1} x_{m-1}^k - a_m x_m^k = 0 \tag{5.4}$$

If the equation (5.4) has at least one natural solution $(x_{10}, ... x_{m0})$, then the number of natural solutions of the equation in $(A)^m$ is equal to:

$$R_m^+(N) = N / sup_{i \leq m}(x_{i0}) = O(N) \tag{5.5}$$

Now we consider the equation:

$$f_2 = a_{m+1} x_{m+1}^k - a_{m+2} x_{m+2}^k + ... + a_{s-1} x_{s-1}^k - a_s x_s^k = 0 \tag{5.6}$$

We assume that: $a_{m+1} \leq ... \leq a_s$. If $(a_{m+2}/a_{m+1})^{1/k}, ..., (a_s/a_{s-1})^{1/k}$ are natural numbers, then based on (4.34) the relation (for the number of natural solutions of the equation (5.6) in $(A)^{s-m}$) is true:

$$R_{s-m}^+(N) \geq N^{(s-m)/2} / (a_{m+2} / a_{m+1}...a_s / a_{s-1})^{1/k} = O(N^{(s-m)/2}). \tag{5.7}$$

If $f_1$ is absent in (5.2), then based on (5.7) we get (for the value $m = 0$):

$$R_s^+(N) \geq O(N^{s/2}). \tag{5.8}$$

If there are $f_1, f_2$, then we obtain (based on assertion 1 ( $f_1 f_2 < 0$) and (5.5)) the following lower bound for the number of natural solutions of the equation (5.2) in $(A)^s$:

$$R_s^+(N) \geq N^{(s-m)/2+1} / sup_{i \leq m}(x_{i0})(a_{m+2} / a_{m+1}...a_s / a_{s-1})^{1/k} = O(N^{(s-m)/2+1}). \tag{5.9}$$

Note that the degree $(s-m)/2+1(m \geq 3)$ is equal to the number of monomials in equation (5.2) with negative coefficients if there are $f_1, f_2$. The number of monomials with negative coefficients in equation (5.2), as mentioned earlier, does not exceed $s/2$. Therefore, the upper bound of the lower estimate for the number of natural solutions of the equation (5.2) is achieved in case then $f_1$ is absence ($m = 0$) and determined by the formula (5.8).

Consider the examples of lower estimates of the number of natural solutions of the equation (5.2).

Example 1 Equation:

$$6x_1^k + 3x_2^k - 3x_3^k - 2x_4^k - 4x_5^k = 0 \tag{5.10}$$

has more than half of monomials with negative coefficients, so we replace the equation (5.10) to the equivalent:

$$-6x_1^k - 3x_2^k + 3x_3^k + 2x_4^k + 4x_5^k = 0 \tag{5.11}$$

We group the monomials of the equation (5.11) as follows:

$$2x_4^k + 4x_5^k - 6x_1^k + 3x_3^k - 3x_2^k = 0. \tag{5.12}$$

Suppose that $f_1 = 2x_4^k + 4x_5^k - 6x_1^k = 0, f_2 = 3x_3^k - 3x_2^k = 0$ in (5.12), then the equation $f_1 = 0$ has the solution $x_4 = x_5 = x_1$ and the number of natural solutions is equal to $R_3^+(N) = N$ in $(A)^3$. Equation $f_2 = 0$ has the solution $x_2 = x_3$ and the number of natural solutions is equal to $R_2^+(N) = N$ in $(A)^2$. Based on assertion 1 ( $f_1 f_2 < 0$) we obtain the following lower estimate of the number of natural solutions of the equation (5.10): $R_s^+(N) \geq N^2 \ln(A)^5$, corresponding to (5.9).

Example 2 Equation:

$$x_1^k + 2x_2^k - 3x_3^k = 0 \tag{5.13}$$

has the form $f_1 = 0$ (function $f_2$ is missing). Indeed, the equation (5.13) has the solution $x_1 = x_2 = x_3$ and the number of natural solutions is equal to: $R_3^+(N) = N$ in $(A)^3$ which corresponds to (5.5).

Example 3 Equation:

$$x_1^k + 2x_2^k - 2x_3 - x_4^k = 0 \tag{5.14}$$

has the form $f_2 = 0$ (function $f_1$ is missing). We group the terms of the equation (5.14) as follows:

$$x_1^k - x_4^k - 2(x_2^k - x_3^k) = 0. \tag{5.15}$$

Equation (5.15) has the solution: $x_1 = x_4, x_2 = x_3$ and the number of natural solutions is equal to: $R_4^+(N) = N^2$ in $(A)^4$ corresponding to (5.8).

Now we consider another class of homogeneous Diophantine equations:

$$c_1 x_1^k + ... + c_s x_s^k = 0, \tag{5.16}$$

where $c_i (2 < i \leq s)$ are integers, $k$ is a natural number and $\Sigma_{i=1}^s c_i = 0$.

In this case, the equation (5.16) can be written as:

$$\Sigma_{1 \leq i < j \leq s} a_{ij}(x_i^k - x_j^k), \tag{5.17}$$

where $a_{ij}$ are non-negative integers.

Equation (5.17) has the solution $x_i = x_j$ if $a_{ij} > 0$. We will estimate the number of natural solutions of the equation (5.17).

In the particular case the equation (5.17) may take the form as:

$$a_{12}(x_1^k - x_2^k) + a_{23}(x_2^k - x_3^k) + \ldots + a_{s-1s}(x_{s-1}^k - x_s^k) = 0, \tag{5.18}$$

which has the solution: $x_1 = x_2 = \ldots = x_{s-1} = x_s$. Based on assertion 1 the equation (5.18) in the area has the following asymptotic lower bound for the number of natural solutions in $(A)^s$:

$$R_s^+(N) \geq O(N). \tag{5.19}$$

In another particular case, the equation (5.17) may take the form as:

$$a_{12}(x_1^k - x_2^k) + a_{34}(x_3^k - x_4^k) + a_{45}(x_4^k - x_5^k) + \ldots + a_{s-1s}(x_{s-1}^k - x_s^k) = 0, \tag{5.20}$$

which has the solution: $x_1 = x_2, x_3 = x_4 = \ldots = x_{s-1} = x_s$. Based on assertion 1, the equation (5.20) has the following asymptotic lower bound of the number of natural solutions in $(A)^s$:

$$R_s^+(N) \geq O(N^2). \tag{5.21}$$

Finally, in the particular case the equation (5.17) may take the form as:

$$a_{12}(x_1^k - x_2^k) + a_{34}(x_3^k - x_4^k) + \ldots + a_{s-1s}(x_{s-1}^k - x_s^k) = 0, \tag{5.22}$$

which has the solution $x_1 = x_2, x_3 = x_4, \ldots, x_{s-1} = x_s$. Based on assertion 1, the equation (5.22) has the following asymptotic lower bound for the number of natural solutions:

$$R_s^+(N) \geq O(N^{[s/2]}), \tag{5.23}$$

where $[B]$ is the integer part of $B$.

Estimate (5.23) is an upper bound of the asymptotic lower estimate for the number of natural solutions of the equation (5.17) in $(A)^s$, which is equal to the upper asymptotic lower estimate for the number of natural solutions of the equation (5.1) in the same area.

Thus, the equation (5.16) (in the case of integers $c_i (2 < i \leq s)$ and $\Sigma_{i=1}^s c_i = 0$) always has solutions in natural numbers with the asymptotic lower estimate of the number of natural solutions in $(A)^s$ changing the value from $O(N)$ to $O(N^{[s/2]})$.

Now let us compare the upper limit of the lower asymptotic estimates for the number of natural solutions of the equation (5.1) with the lower bound of the upper asymptotic estimate of the number of natural solutions of the equation obtained by the CM.

It is known that the CM is just to get the upper estimates for the number of natural solutions of the equation (5.1) at $s > g(k)$. Based on Vinogradov's theorem the value $s > 2^k$ is just for $k \geq 8$ and the function $g(k)$ increases faster than the exponential function at lower values $k$.

Based on (3.5) the upper asymptotic estimate of the number of natural solutions of the equation (5.1), obtained by the CM in $(A)^s$ is:

$$O(N^{[s/2]}) \leq R_s^+(N) << N^{[s/2]+\epsilon}, \tag{5.24}$$

where $\epsilon$ is a small real number.

If $s > 2^k$ then it is true $k \leq [\log_2(s)] < [s/2]$. Consequently, it is performed $k < [s/2]$ in (5.24). Therefore, the lower limit of the upper asymptotic estimate for the number of natural solutions of the equation (5.1) in $(A)^s$ is:

$$R_s^+(N) << N^{s-k+\epsilon} < N^{[s/2]+\epsilon}. \tag{5.25}$$

Since the upper limit of the lower asymptotic estimate of the number of natural solutions of homogeneous Diophantine diagonal equations (4.41), (5.1), (5.16) in $(A)^s$: $R_s^+(N) \geq O(N^{[s/2]})$ and using (5.25) we obtain the following limits of estimates:

$$O(N^{[s/2]}) \leq R_s^+(N) << N^{[s/2]+\epsilon}. \tag{5.26}$$

## 6. CONCLUSION AND SUGGESTIONS FOR FURTHER WORK

The next article will focus on the estimations of the number of solutions of inhomogeneous algebraic Diophantine equations with integer coefficients.

## 7. ACKNOWLEDGEMENTS

Thanks to everyone who has contributed to the discussion of this paper.

# References

1. Stepanov SA , The arithmetic of algebraic curves, 1991, 185 p.

2. Faltings G. Endlishkeitssätzt für abenlsche Varietäten über Zahlköpern. Invent. Math.ß 1983.V.73, №3.-P. 349-366.

3. Siegel C.L. Approximation algebraisher Zahlen. Math. Zeitschr.- 1921.-Bd 10. –S. 173-213

4. Schmidt W.M. Norm form equations. Ann. Math. - 1972.-V. 96. - P. 526-551.

5. Baker A. Linear form in the logarithms of algebraic numbers. Mathematika. – 1966. –V. 13 – P. 204-216; 1967. - V. 14. P. 102-107, 220-228; 1968. – V.15. – P. 204-216.

6. Sprindzhuk VG Classic Diophantine equations in two unknowns. - M .: Nauka, 1982.

7. Baker A. Bounds for the salutations of the hyperelliptic equation. Proc. Cambr. Philos. Soc. -1969. – V. 65. – P. 439-444

8. V.G. Sprindzuk Hyperelliptic Diophantine equations and class number. Acta Arithm. - 1976. V. 30, №1. - P. 95-108

9. R. Vaughan, "The method of Hardy and Littlewood", M, Mir, 1988, 184 p.

10. I.M. Vinogradov The method of trigonometric sums in number theory. M .: Nauka, 1972.

11. Victor Volfson. Estimations of the number of solutions of algebraic Diophantine equations with natural coefficients using the circle method of Hardy-Littlewood. arXive preprint http://arxiv.org/abs/1512.07079 (2015)

12. V. Sierpinski "The solution of equations in integers" Publishing house "Fizmatgiz", M., 1961, 88 p.

13. Buchstab "Number Theory", Publishing House of the "Enlightenment", M, 1966, 384 p.

14. D. R. Heath-Brown. The density of rational points on nonsingular hypersurfaces. Proc. Indian Acad. Sci. Math. Sci., 104(1):13–29, 1994.

15. D. R. Heath-Brown. Sums and differences of three kth powers. J. Number Theory, 129(6):1579–1594, 2009.